\documentclass[11pt,amssym,twoside]{article}
\usepackage{amssymb}
\newtheorem{thm}{Theorem}[section]
\newtheorem{mthm}[thm]{Main Theorem}
\newtheorem{cor}[thm]{Corollary}

\newtheorem{conj}[thm]{Conjecture}
\newtheorem{prop}[thm]{Proposition}
\newtheorem{rthm}[thm]{Riemann's Existence Theorem}
\newtheorem{defn}[thm]{Definition}

\newtheorem{rem}[thm]{Remark}

\newcommand{\To}{\longrightarrow}

\newcommand{\Z}{\mathbb Z}

\begin{document}

\title{On Profinite Hyperbolicity and Diophantine Geometry}
\author{Arash Rastegar}


\maketitle
\begin{abstract}
In this note, we explore the notion of hyperbolicity of topologically finitely
generated profinite groups. Some
applications to diophantine geometry are suggested and we try to 
reformulate certain problems in diophantine
geometry in terms of hyperbolic profinite groups.
Then, we introduce many occasions in which Galois groups
are free profinite and try to explore implications of this condition 
in the world of diophantine geometry. In particular, we prove that,
Grothendieck's "section conjecture" plus Shafarevich's "freeness
conjecture" imply that hyperbolic curves have infinitely many
solutions over the maximal abelian extension of a global field.
This makes Moerdell's conjecture, which was proved by Faltings,
more interesting.
\end{abstract}

\section*{Introduction}
Gromov's theory of hyperbolic groups, studies finitely generated
groups whose Cayley graphs are hyperbolic metric spaces, when
equipped with word metric. This notion is independent of choices
made for finitely many generators. We intend to explore the implications of extension of this theory
to finitely generated profinite groups. The
purpose of this paper is to imply that hyperbolicy of profinite
groups has applications in finiteness and infiniteness results in diophantine
geometry.

Given an elliptic curve $E$ over $\mathbb{Q}$ and fixing a prime
$\ell$, Tate  associated a representation
$$
\rho_{\ell}:Gal(\overline{\mathbb{Q}}/\mathbb{Q})\rightarrow
GL_{2}(\mathbb{Z}_{\ell})
$$
by taking the inverse limit over $n$ of the groups of
$\ell^{n}$-torsion points of $E$ as Galois modules. On the other
hand, one can think of the $\ell^{n}$-torsion points of $E$ as
the Galois group of the function field extension associated to
the multiplication map $\times\ell^{n}:E \rightarrow E$. Since
every Galois etale covering of an elliptic curve is again an
elliptic curve, we can take a limit over all Galois etale covers
$E'\rightarrow E$ of degree $\ell^{n}$ for some $n$, and obtain
the same Galois module structure on $\mathbb{Z}_{\ell}\times
\mathbb{Z}_{\ell}$ as part of the algebraic fundamental group of
$E$.

For an arbitrary curve $X$ defined over an algebraically closed field of characteristic zero,
the algebraic fundamental group
$\pi_{1}^{alg}(X)$ is defined by Grothendieck as $\lim Gal(K'/K)$
where $K$ is the function field of $X$ and $K'$ runs over all
Galois extensions of $K$ such that the corresponding curve $X'$
is etale over $X$. For example, $\pi^{alg}_{1}(\mathbb{P}^{1})=1$
and for every elliptic curve $E$ we have
$$
\pi^{alg}_{1}(E)=\prod_{\ell} \mathbb{Z}_{\ell}\times \mathbb{Z}_{\ell}.
$$

Grothendieck proved that, for a curve $X$ of genus $g$, the
algebraic fundamental group of $X$ is isomorphic to the
completion of the ordinary topological fundamental group of $X$
over $\mathbb{C}$ [Grt]. If $X$ is defined over $\mathbb{Q}$ one can
induce an action of $Gal(\overline{\mathbb{Q}}/\mathbb{Q})$ on
$\pi^{alg}_{1}(X)$.

One can also define the $\ell$-adic algebraic fundamental group
$\pi^{alg}_{1,\ell}(X)$ as the limit
$$
\lim Gal(K'/K)
$$
where $K$ is the function field of $X$ and $K'$ runs over all
degree $\ell^{n}$ Galois extensions of $K$ such that the
corresponding curve $X'$ is etale over $X$. The group
$\pi^{alg}_{1,\ell}(X)$ is isomorphic to the $\ell$-adic
completion of $\pi^{top}_{1}(X)$. If $X$ is defined over
$\mathbb{Q}$ one can associate an $\ell$-adic Galois
representation to $X$ which is a direct summand of the above
mentioned representation
$$
\rho_{\ell}:Gal(\overline{\mathbb{Q}}/\mathbb{Q})\rightarrow
Aut(\pi^{alg}_{1,\ell}(X)).
$$

The Frattini subgroup $\Phi(G)$ of the pro-$\ell$ group
$G=\pi^{alg}_{1,\ell}(X)$ is fixed by each automorphism [Rib-Zal].
Therefore we get a representation
$$
\rho_{\ell}:Gal(\overline{\mathbb{Q}}/\mathbb{Q})\rightarrow
Aut(G/\Phi(G)).
$$
Since $\Phi(G)=\overline{G^{\ell}[G,G]}$ we can think of $G/\Phi(G)$
as a finite dimensional vector space over the field
$\mathbb{F}_{\ell}$ with $\ell$ elements [Rib-Zal]. This way we have
recovered the mod-$\ell$ Galois representation associated to $X$ via
Tate module.

The classical approach to definition of an $\ell$-adic Galois
representation is to consider the action of Galois group on the
Tate module associated to the Jacobian variety $Jac(X)$ as above.
This representation can be reconstructed from the one we
introduced above by abelianization. But lots of
information is lost during this process. We believe the structure of
fundamental group can help us to understand the Galois action
better.

The notion of etale fundamental group introduced by Grothendieck,
links the notion of Galois group with its function field
counterpart, which is geometric fundamental group [Grt]. By
"Grothendieck's anabelian conjecture" the arithmetic fundamental
group of a hyperbolic algebraic curve completely determines the
algebraic structure of the curve. More precisely, if $X$ is an
algebraic variety over a field $k$  which is proper, connected, and geometrically reduced,
we have the following exact
sequence
$$
0 \To \pi_1(X_{\bar {k}}) \To \pi_1(X_k) \To Gal(\bar{k}/k) \To 0
$$
where $\pi_1(X_{\bar {k}})$ is called the geometric fundamental group and
$\pi_1(X_k)$ is called the arithmetic fundamental group. This determines
a homomorphism
$$
\pi_1(X_k) \to Aut(\pi_1(X_{\bar {k}}))
$$
by conjugating the group $\pi_1(X_{\bar {k}})$ by elements of
$\pi_1(X_k)$ and then, after taking quotients, one obtains a
homomorphism landing in the outer automorphism group
$$
\rho_X: Gal(\bar{k}/k) \to Out(\pi_1(X_{\bar {k}}))
$$
To consider the outer action of Galois group is equivalent to
considering $\pi_1(X_k)$ as an extension of $Gal(\bar{k}/k)$.

By Grothendieck's "fundamental conjecture", if $k$ is finitely
generated over the prime field and if $X$ is a hyperbolic
algebraic variety, or in other words sufficiently non-abelian,
$\rho_X$ determines $X$ over $k$. Grothendieck's "hom conjecture"
states that if $X$ and $Y$ are algebraic varieties over a finitely
generated field $k$, then the natural map between dominant
morphisms from $X$ to $Y$ over $k$ to Galois compatible open
homomorphisms between the arithmetic fundamental groups
$$
Hom(X,Y) \to Hom_{Gal(\bar{k}/k)}(\pi_1(X_k),\pi_1(Y_k))
$$
induces a bijection if right hand side is considered up to
composition with an inner automorphism of $\pi_1(Y_k)$. According
to "section conjecture", splittings of the previously mentioned short exact
sequence are in one to one correspondence with the $k$-rational
points of the projectivized variety. Moreover, sections
corresponding to points at infinity are group theoretically
characterized. The are several more explicit versions of this
general philosophy, similar to "fundamental conjecture", "hom
conjecture" and the "section conjecture". We refer the reader to
[Nak-Tam-Moc] for an exposition.

The known half of the section conjecture, tells us that to any
point of $X$ defined over $k$ one can associate a splitting 
and this association is injective. Therefore, in order to
prove finiteness of the rational points it is enough to show that
$\pi_1(X_k)$ has only finitely many splittings as an extension of
$Gal(\bar{k}/k)$ by $\pi_1(X_{\bar {k}})$. This brings us to the
realm of profinite group theory. We suggest that hyperbolicity of
$\pi_1(X_{\bar {k}})$ is enough to get such finiteness results.
This leads us to the following conjecture which implies Mordell's
conjecture proved by Faltings:

\begin{conj}
Let $H$ denote a finitely generated profinite hyperbolic group,
and let $G$ denote a split extension of the Galois group of 
the maximal extension unramified outside a set $S$ of finitely many primes
of a number field $k$ by a hyperbolic finitely generated profinite group $H$
$$
0 \to H \to G \to Gal(\bar{k}_S/k) \to 0
$$
Then, there are only finitely many $H$-conjugacy classes of
splittings of the above short exact sequence.
\end{conj}

In fact, we expect that $Gal(\bar{k}_S/k)$ be profinite hyperbolic and the conjecture
to hold for Galois group replaced by any profinite hyperbolic group. 
Note that, we have the following finiteness result for restricted
ramification (10.9.11 in [Neu-Sch-Wij]).
For $k$ a number field and $S$ a finite set of primes in
$K$, then $Gal_S(K)$ is topologically generated by conjugacy
classes of finitely many elements.

We also try to imply that freeness of Galois groups has
important implications in diophantine geometry. We will use the
freeness of Galois groups to get infiniteness results. The
following conjecture of Shafarevich states that certain Galois
groups are profinite free:

\begin{conj}(Shafarevich)
Let $k$ be a global field and $k^{cycl}$ denote its maximal
abelian extension. Then $Gal(k^{cycl}/k)$ is a free profinite
group.
\end{conj}
There are a number of evidences for Shafarevich's
conjecture. Iwasawa proved that Galois group of the maximal
solvable extension of $k^{cycl}$ is profinite free, and Tate
proved that $Gal(k^{cycl}/k)$ has cohomological dimension one
[Ser] hence by a theorem of Gruenberg it is projective [Gru].

Shafarevich's conjecture was proved in the function field case by F. Pop [Pop1]. We
shall see that Grothendieck's "anabelian" considerations plus
Shafarevich's conjecture imply infiniteness of points defined
over the maximal unramified extensions of numer fields which
brings Mordell's conjecture proved by Falthings to a new
perspective.

\begin{mthm}
Let $k$ denote a global field and $X$ a hyperbolic algebraic curve
defined over $k$, which has a point over $k$. Then, Grothendieck's "section conjecture" plus
Shafarevich's "freeness conjecture" imply that the number of
points of $X$ defined over $k^{cycl}$ is infinite.
\end{mthm}

\begin{rem}
 If Shafarevich's conjecture is true for any finitely generated field in characteristic zero
we get the same result over such fields.
\end{rem}

Organization of the paper is as follows. In the first section, we
review basic definition of hyperbolic metric spaces and
hyperbolic groups. We review
needed background in profinite groups and pro-$p$ groups as a
special class of profinite groups in the second section. In the
third section, we define finitely generated hyperbolic profinite
groups. The forth section contains a profinite hyperbolic group
theory formulation of diophantine geometry.
in the fifth section, we review review results on topological fundamental group
and in the sixth section, we introduce the algebraic geometric formulation of the
Poincare conjecture. Then we review the construction of outer representations. 
In next section, we show that outer representation lands in the completion of the mapping class group.
In the ninth section, we
review basic definition of free profinite groups and profinite
free products. Then, in the tenth section, we review a few
classical results in algebraic number theory which imply that
certain Galois groups are profinite free. Afterwards, we explore
connections with diophantine geometry and prove the main theorem
in the final section.

\section{Preliminaries on hyperbolic groups}

Motivated by hyperbolic phenomena in the geometry of manifolds
and spaces of negative curvature, Gromov defined the notion of
hyperbolic finitely generated groups. Trying to make the
formulation of this concept independent of curvature, he invented
the notion of a quasi-isometry between metric spaces [Grm].

\begin{defn}
Let $(X,d)$ be a metric space. The Gromov product of points $x$
and $y$ of $X$ with respect to a base point $x_0$ is defined to
be the real number
$$
(x.y)=(x.y)_{x_0}=\frac{1}{2}(d(x,x_0)+d(y,x_0)-d(x,y))
$$
The metric space $X$ is called `Gromov hyperbolic' or simply
`hyperbolic', if it is $\delta$-hyperbolic for some $\delta \geq
0$ which means that
$$
(x.y)\geq min((x.z),(y.z))-\delta
$$
should hold for all $x,y,z\in X$ and every choice of basepoint
$x_0$. A geodesic metric space $X$ is hyperbolic if and only if
for some $\delta \geq 0$ every geodesic triangle is $\delta$-thin,
meaning that every side is contained in the $\delta$-neighborhood
of the union of two other sides.
\end{defn}

Important examples of hyperbolic spaces are: the hyperbolic upper
half-space $\Bbb H^n$ and all simply-connected Riemannian
manifolds with sectional curvature bounded above but away from
zero, $\Bbb R$-trees and Cayley graphs of hyperbolic groups.

\begin{defn}
A finitely generated group $G$ is called hyperbolic, if for some
finite generating set $\Gamma$ the Cayley graph of $G$ with
respect to $\Gamma$ is a hyperbolic space with word metric.
\end{defn}

In order to show that hyperbolicity is well-defined one has to use
the notion of quasi-isometries between metric spaces and the fact
that hyperbolicity is a quasi-isometric invariant. Important
examples of hyperbolic groups are finite groups, free groups of
finite rank, fundamental groups of compact Riemannian manifolds
with negative sectional curvature, and groups properly
discontinuously acting on hyperbolic spaces with compact
quotient. Here are some properties of hyperbolic groups: A
hyperbolic group is finitely presented. It can not contain a
subgroup isomorphic to $\Bbb Z\oplus \Bbb Z$. A finite index
subgroup is hyperbolic if and only if the ambient group is. A
hyperbolic groups acts freely and transitively on a contractible
finite dimensional locally finite polyhedral complex. As a
consequence, it contains only finitely many conjugacy classes of
torsion elements.

\section{Preliminaries on profinite groups}

Interest in profinite groups originates in the study of the
Galois groups of infinite extensions of fields. Galois groups
naturally carry `Krull' topology under which they become
Hausdorff compact totally disconnected topological groups. These
properties characterize profinite groups and indeed, profinite
groups are precisely Galois groups.

\begin{defn}
Let $C$ denote a class of finite groups. A pro-$C$ group $G$ is an
inverse limit of a surjective inverse system of groups $G_i$ in
$C$ with discrete topology.
$$
G=\lim_{\leftarrow} G_i
$$
We say that $G$ is topologically generated by a finite subset
$\Gamma$ if its image in $G_i$ generates $G_i$, or equivalently if
the closure of the group generated by $\Gamma$ is $G$. In this
case, we say that $G$ is the pro-$C$ completion of the group
generated by $\Gamma$.
\end{defn}

In case $C$ is the class of all finite groups, then $G$ is called
a profinite group and if $C$ is the class of $p$-groups, then $G$
is called a pro-$p$ group. Since we intend to discuss hyperbolic
profinite groups, we shall first recognize the class of free
pro-$C$ groups as the most important class of hyperbolic
profinite groups, to be introduced in the next section.
\begin{defn}

For any profinite space $X$ there exists a unique pro-$C$ group
$F_C(X)$ containing $X$ as a profinite subset which satisfies the
following universal property: Any continuous map $X \to G$ to a
pro-$C$ group $G$ whose image generates $G$, uniquely breaks
through inclusion $X\hookrightarrow F$ and a continuous
homomorphism $F\to G$. $F_C(X)$ is called the free pro-$C$ group
on the profinite basis $X$.
\end{defn}

In the case of finitely generated free pro-$C$ groups, any
generating set with the same cardinality as the basis is again a
basis. In fact, pro-$C$ completion of the free group generated by
finite basis $X$ is the free pro-$C$ group on the basis $X$. Free
pro-$C$ groups are examples of free pro-$C$ products.

\begin{defn}
A family of homomorphisms $\rho_i :H_i\to H$ between pro-$C$
groups, for $i\in I$, is called convergent to $1$ if $\rho_i
(H_i)$ is convergent to $1$ which means that any open
neighborhood of $1$ contains almost all $\rho_i (H_i)$. The free
pro-$C$ product of a family of $H_i$ is a pro-$C$ group
$H=*_{i\in I} H_i$ together with a convergent family of
homomorphisms $\rho_i :H_i\to H$ which is universal among all
such families targeting in any group $H'$. Free pro-$C$ products
are unique up to isomorphisms.
\end{defn}

Subgroups of finite index in free pro-$C$ products are again free
pro-$C$ products, which could be an indication of the fact that
free products being similar to free groups are usually hyperbolic
groups. For example, free products of finitely many finitely
generated hyperbolic groups are again hyperbolic groups.

\section{Hyperbolic profinite groups}

Now, we shall try to introduce a profinite version of the
hyperbolic phenomena. In order to discuss hyperbolicity of a
profinite group, we need to compensate existence of a geometric object which replaces the role
of Cayley graph of a hyperbolic group.
\begin{defn}
Let $G$ be a residually finite profinite group, which is
isomorphic to inverse limit of all its finite images $G_i$ of
$G$, and let $\Gamma$ denote a finite topological generating set. 
If $\Gamma$ generates a  hyperbolic finitely generated subgroup, we say that $G$ is
profinite hyperbolic.
\end{defn}

\begin{prop}
Free pro-$C$ products of finitely many finitely hyperbolic pro-$C$
groups $G_i$ are again hyperbolic pro-$C$ groups.
\end{prop}
\textbf{Proof.} Let $\Gamma_i$ denote a finite generating set for
$G_i$ for $i\in I$. If the Cayley graphs of $<\Gamma_i>\subset G_i$ with respect to
$\Gamma_i$ are $\delta_i$-hyperbolic, then the Cayley graph of
$*_{i\in I} G_i$ with respect to $\cup_i\Gamma_i$ will also be
$\delta$-hyperbolic for some $\delta$. $\Box$

\section{Profinite hyperbolicity and Diophantine objects}

We shall start with the following basic conjecture:
\begin{conj}
For any number field $k$, the Galois group $Gal(\bar{k}_S/k)$ of 
the maximal extension unramified over a set $S$ of finitely many primes
of a number field $k$ is a hyperbolic profinite group.
\end{conj}

There exists a pro-$p$ version of the Grothendieck short exact
sequence for a curve defined over any finitely generated field
$$
0 \To \pi^{(p)}_1(X_{\bar k}) \To \pi^{(p)}_1(X_k) \To
Gal(\bar{k}/k)(p) \To 0
$$
and the philosophy of Grothendieck's section conjecture could be
also implied in this case. In a forthcoming section, many examples
are introduced where we get hyperbolicity of certain pro-$p$
Galois groups for $p$ outside a finite set of primes $S$. Together
with hyperbolicity of the corresponding $\pi^{(p)}_1(X_{\bar k})$
one is lead to the following conjecture

\begin{conj}
Let $k$ denote a finitely generated field and $X$ a hyperbolic
algebraic curve defined over $k$, and assume that
$Gal(\bar{k}_S/k)(p)$ is hyperbolic for the finite set
of primes $S$ of bad reduction for $X$, then the number of points of $X$ defined over
$A_S$ is infinite, where $A_S$ denotes the integral closure of $\Z$
in the maximal extension of $k$ unramified outside $S$.
\end{conj}

\begin{rem}
Let $k$ denote a finitely generated field and $X$ a hyperbolic
algebraic curve defined over $k$ of genus $\geq 1$ after compactification.
We prove prove in [Ras] that the number of points of $X$ defined over
$A_S$ is infinite, where $A_S$ denotes the integral closure of $\Z$
in the maximal extension of $k$ unramified outside a finite set of places $S$.
\end{rem}

In a more abstract setting, one can also think of a profinite
group $\Gamma$ surjecting to the Galois group, as an "arithmetic"
space
$$
\Gamma \To Gal(\bar{k}/k) \To 0
$$
and ask for pure algebraic hyperbolicity conditions implying the
finiteness of splittings $Gal(\bar{K}/K) \to \Gamma$ of the
corresponding short exact sequence. This would be a
generalization of the Mordell conjecture. Using such algebraic
methods one will have a proof of the following folklore more
strong version of Mordell's conjecture:

\begin{conj}
The number of rational points of any hyperbolic curve defined
over a number field, is bounded above by a bound depending only on
the genus and the number field.
\end{conj}

This result is proved in the function field case using Model
theory. Up to know there exists no number theoretical proof of the
Mordell's conjecture in the function field case. This is one of the main implications
of the section conjecture.

\section{The topological fundamental group revisited}

The fundamental group plays a fundamental role in geometry. 
Let $X$ be a compact surface of genus $g$.
Hempel in 1976 proved that Poincar\'e conjecture holds iff there exists a
unique epimorphism
$$
\pi^{top}_{1}(X)\rightarrow F(g)\times F(g)
$$
up to automorphisms of $\pi^{top}_{1}(X)$ and automorphisms of
the components $F(g)$ which are free groups with $g$ generators
[Hem]. Existence of such an epimorphism is evident from the
presentation
$$
\pi^{top}_{1}(X)=\langle
a_{1},...,a_{g},b_{1},...,b_{g}|\prod_{i=1}^{g}[a_{i},b_{i}]=1
\rangle.
$$

Having this presentation, one can show that the fundamental group
is a torsion-free group with trivial center. We also have access
on a necessary and sufficient condition for a one-relator group to
be isomorphic to the fundamental group [Col-Zie].

For a free group $F$, we say that two surjective maps
$\varphi_1,\varphi_2:\pi^{top}_{1}(X)\rightarrow F$ are
equivalent if there exist an automorphism of $\pi^{top}_{1}(X)$
and an automorphism of $F$ which take $\varphi_1$ to $\varphi_2$,
i.e. there exists $\alpha\in Aut(\pi^{top}_{1}(X))$ and $\beta\in
Aut(F)$ such that the following diagram is commutative
$$
\begin{array}{ccc}
\varphi_1:\pi^{top}_{1}(X)& \longrightarrow  & F \\
\alpha\downarrow &  & \downarrow\beta \\
\varphi_2:\pi^{top}_{1}(X)&\longrightarrow  & F \\
\end{array}
$$
We say that $\varphi_1$ and $\varphi_2$ are strictly equivalent
if there exist an automorphism $\alpha$ of $\pi^{top}_{1}(X)$
such that $\varphi_2=\varphi_1\circ\alpha$. Zieschang in 1964
proved that if rank of $F$ is less than or equal to $g$ there
exists only one equivalence class of surjective maps
$\pi^{top}_{1}(X)\rightarrow F$ [Zie]. With the same assumption,
Kurchanov and Grigorchuk in 1989 proved that there exists only
one strict equivalence class of such surjective maps [Kur-Gri].
Existence of such a surjection is an easy consequence of the
presentation of the fundamental group. Nielsen in 1927 proved
that every automorphism of $\pi^{top}_{1}(X)$ is induced by an
automorphism of the corresponding free group which preserves the
relation $R$ [Nie].

Suppose $H_0$ and $H_1$ are subgroups of $\langle
a_{1},...,a_{g},b_{1},...,b_{g}|\prod_{i=1}^{g}[a_{i},b_{i}]=1
\rangle$ which are generated by $a_1,...,a_g$ and $b_1,..,b_g$
respectively. The subgroups $L^+$ and $L^-$ defined as such
$$
L^-=\{\varphi\in Aut(\pi^{top}_{1}(X)):\varphi(H_0)=H_0\}
$$
$$
L^+=\{\varphi\in Aut(\pi^{top}_{1}(X)):\varphi(H_1)=H_1\}
$$
are conjugate in $Aut(\pi^{top}_{1}(X))$. The group
$Aut(\pi^{top}_{1}(X))$ is generated by $L^-$ and $L^+$.

Every automorphism of $\pi^{top}_{1}(X)$ induces an automorphism
of the abelian group
$\pi^{top}_{1}(X)/[\pi^{top}_{1}(X),\pi^{top}_{1}(X)]$ whose
matrix coordinatizes a representation
$$
\rho:Aut(\pi^{top}_{1}(X))\longrightarrow Sp_{2g}(\mathbb{Z})
$$
Here $Sp_{2g}(\mathbb{Z})$ is the group of matrices $S$ of order
$2g$ such that $SJS^t=\pm J$ where
$$
J=\left(\begin{array}{cc}
0 & I_g \\ -I_g & 0
\end{array} \right).
$$
The representation $\rho$ takes elements of $L^-$ to lower
triangular matrices
$$
\left(\begin{array}{cc} a & 0 \\ b & c
\end{array} \right)
$$
and elements of $L^+$ go to upper triangular matrices
$$
\left(\begin{array}{cc} a & b \\ 0 & c
\end{array} \right).
$$
\section{The algebro-geometric formulation of Poincar\'e conjecture}

Providing information about the completed fundamental group is the
key to a better understanding of the Galois representation. Let
$G$ denote the pro-$p$ completion of the fundamental group. Taking
pro-$p$ completions from the surjective maps
$\pi^{top}_{1}(X)\rightarrow F(g)\times F(g)$ and
$\pi^{top}_{1}(X)\rightarrow F(2g)$ one gets surjections from $G$
to the pro-$p$ completions $F(g)_{p}\times F(g)_{p}$ and
$F(2g)_{p}$. Also, every automorphism of $\pi^{top}_{1}(X)$
induces an automorphism of $G$. One can ask if the analogue of
the Poincar\'e conjecture is true here. 
\begin{conj}
There is only one strict equivalence class of surjections
$$
G \longrightarrow F(g)_{p}\times F(g)_{p}.
$$
\end{conj}

The representation $\rho:Aut(\pi^{top}_{1}(X))\longrightarrow
Sp_{2g}(\mathbb{Z}) $ induces a map
$$
\widehat{\rho}:Aut(G) \longrightarrow Sp_{2g}(\mathbb{Z})_p
$$
which after combining with the natural surjection
$$
Sp_{2g}(\mathbb{Z})_p\longrightarrow Sp_{2g}(\mathbb{Z}_p)
$$
induces a representation
$$
\rho_X:Aut(G) \longrightarrow Sp_{2g}(\mathbb{Z}_p).
$$
If we combine $\rho$ with the Galois action on $\pi^{alg}_{1}(X)$
one has associated a $p$-adic Galois representation to the
algebraic curve $X$:
$$
\rho_p:Gal(\overline{\mathbb{Q}}/\mathbb{Q})\rightarrow
Sp_{2g}(\mathbb{Z}_p).
$$
This is the same as the Galois representation associated to the
Tate module of the Jacobian variety of $X$.

\section{The Galois action on $\pi^{alg}_{1}(X)$}

The outer automorphism group of $\pi^{alg}_{1}(X)$ as a quotient
of $Aut(\pi^{alg}_{1}(X))$ accepts a representation from the
Galois group
$$
\rho_p:Gal(\overline{\mathbb{Q}}/\mathbb{Q})\longrightarrow
Aut(\pi^{alg}_{1}(X))\longrightarrow Out(\pi^{alg}_{1}(X))
$$
Mochizuki proved that if $X$ and $X'$ are curves over
$\mathbb{Q}$ of genus greater than or equal to $2$, the map
$$
Isom_{\mathbb{Q}}(X,X')\longrightarrow
Out_{Gal(\overline{\mathbb{Q}}/\mathbb{Q})}
(Out(\pi^{alg}_{1}(X)),Out(\pi^{alg}_{1}(X')))
$$
is a one-to-one correspondence, where
$Out_{Gal(\overline{\mathbb{Q}}/\mathbb{Q})}$ denotes Galois
equivariant isomorphisms between the two groups [Moc]. In
particular, the Galois representation on $Out(\pi^{alg}_{1}(X))$
completely determines the curve $X$ over $\mathbb{Q}$.

The representation $\widehat{\rho}$ breaks through
$Out(\pi^{alg}_{1}(X))$ because inner automorphisms are neutral
when reduced to abelianization of a group
$$
\widehat{\rho}:Aut(\pi^{alg}_{1}(X)) \longrightarrow
Out(\pi^{alg}_{1}(X))\longrightarrow Sp_{2g}(\mathbb{Z})_p.
$$
Thus for each prime $p$ one can associate a $p$-adic
representation to $Out(\pi^{alg}_{1}(X))$ by considering the
pro-$p$ part of the outer automorphism group
$$
\rho_X:Out(\pi^{alg}_{1}(X)) \longrightarrow
Sp_{2g}(\mathbb{Z}_p).
$$
and therefore the Galois representation
$$
\rho_p:Gal(\overline{\mathbb{Q}}/\mathbb{Q})\rightarrow
Sp_{2g}(\mathbb{Z}_p)
$$
could be associated to the Galois action on
$Out(\pi^{alg}_{1}(X))$. This shows that the language of outer
automorphisms is an appropriate one for studying the arithmetic
of curves over $\mathbb{Q}$. To study the geometric Galois
representation
$$
\rho:Gal(\overline{\mathbb{Q}}/\mathbb{Q})\rightarrow
Out(\pi^{alg}_{1}(X))
$$
the first step, is to study $Out(\pi^{alg}_{1}(X))$.

\section{Preliminaries on free profinite groups}

We shall first define the class of free pro-$C$ groups via its
universal property:

\begin{defn}
For any profinite space $X$ there exists a unique pro-$C$ group
$F_C(X)$ containing $X$ as a profinite subset which satisfies the
following universal property: Any continuous map $X \to G$ to a
pro-$C$ group $G$ whose image generates $G$, uniquely breaks
through inclusion $X\hookrightarrow F$ and a continuous
homomorphism $F\to G$. $F_C(X)$ is called the free pro-$C$ group
on the profinite basis $X$.
\end{defn}

In the case of finitely generated free pro-$C$ groups, any
generating set with the same cardinality as the basis is again a
basis. In fact, pro-$C$ completion of the free group generated by
finite basis $X$ is the free pro-$C$ group on the basis $X$. Free
pro-$C$ groups are examples of free pro-$C$ products.

One can obtain free profinite groups by taking profinite free
product of free profinite groups of rank one. This is because
profinite free products of profinite free groups are again
profinite free. Here we remind the notion of profinite free
product.

\begin{defn}
A family of homomorphisms $\rho_i :H_i\to H$ between pro-$C$
groups, for $i\in I$, is called convergent to $1$ if $\rho_i
(H_i)$ is convergent to $1$ which means that any open
neighborhood of $1$ contains almost all $\rho_i (H_i)$. The free
pro-$C$ product of a family of $H_i$ is a pro-$C$ group
$H=*_{i\in I} H_i$ together with a convergent family of
homomorphisms $\rho_i :H_i\to H$ which is universal among all
such families targeting in any group $H'$. Free pro-$C$ products
are unique up to isomorphisms.
\end{defn}

\section{Free Galois groups}

This section contains no original results. All of the theorems
mentioned are classical results of algebraic number theory. We
have just collected a few examples of profinite free Galois groups
to be used in the next section.

Let $K(C)$ denote the maximal pro-$C$ extension of $K$ and
$Gal_K(C)$ denote the maximal pro-$C$ quotient of $Gal_K$. For a
global field $K$ and a prime $p$, considering a class of
$p$-groups, we have an injection $Gal_{K_{\beta}}(p) \to
Gal_K(p)$ where $\beta$ is a prime of $K(p)$. If
$\beta_1,...,\beta_r$ are different primes of $K(p)$ then the
natural map $*G_{K_{\beta_i}}(p) \to G_K(p)$ is injective
(Theorem 9.3.2 in [Neu-Sch-Wij]). This could imply freeness in
certain cases. Let us discuss the function field case first. \\
$\bold 1.$ Let $K$ be a function field of characteristic $p$ and
let $k\subset K$ be the finite constant field with algebraic
closure $\bar k$. Let $S$ be a finite set of places of $K$ of
cardinality $n$ and $Gal_S(K)=Gal(K_S/K)$. We have the following
exact sequence:
$$
0\to Gal_S(K\bar{k})\to Gal_S(K)\to Gal_k\to 0
$$
$Gal_k$ is a free profinite group of rank $1$ topologically
generated by Frobenius automorphism. The group $Gal_S(K\bar{k})$
is trivial in case $S=\emptyset$ and zero genus $g(K)=0$, and thus
$Gal_S(K)= Gal_k\cong \hat{\Bbb Z}$. In case $S\neq\emptyset$,
for every $l\neq p$ we have a pro-$l$ presentation of
$Gal_S(K\bar{k})(l)$ by $2g+n$ generators and one relation
$$
Gal_S(K\bar{k})(l)\cong \langle x_1,...,x_{2g},y_1,...,y_n|
[x_1,x_2]...[x_{2g-1},x_{2g}]y_1...y_n=1 \rangle.
$$
which is a free pro-$l$ group of rank $2g+n-1$. In fact the
following is true:
\begin{rthm}
(Algebraic Form) Let $K$ be function field of characteristic $p$
with finite constant field $k$ and let $S$ be a non-empty set of
places of $K$ and $g(K)>0$ then for every $l$ different from $p$
we have the following exact sequence
$$
0\to *_{p\in S(K_{\emptyset}(p))} T_p(l)\to Gal_S(K\bar{k})(l)\to
G_{\emptyset}(K\bar{k})(l)\to 0
$$
where $*$ denotes the free pro-$l$ product and $T_p\cong \Bbb
Z_l(1)$ is the pro-$l$ quotient of inersia at $p$ and
$S(K_{\emptyset}(p))$ denotes the set of primes above $S$ in
$K_{\emptyset}(p)$.
\end{rthm}
$\bold 2.$ A number theoretical analogue of Riemann's existence
theorem exists:
\begin{rthm}
(Number Theoretical Form) Let $K$ be number field and $p$ be a
prime, and let $S\supset T \supset S_p(K)\cup S_{\infty}$ be sets
of places of $K$, then the canonical homomorphism
$$
 *_{p\in S\setminus T(K_T(p))} T(K_p(p)/K_p)\to Gal(K_S(p)/K_T(p))
$$
is an isomorphism where $*$ denotes the free pro-$l$ product and
$T(K_p(p)/K_p)\subset Gal(K_p(p)/K_p)$ is the inertia group at $p$
and $S\setminus T (K_T(p))$ denotes the set of primes above
$S\setminus T$ in $K_T(p)$.
\end{rthm}

This is again a freeness assertion. Note that, for primes $p\in
S\setminus T(K)$ for which $N(p)\equiv 1$ mod $p$ we have
$T(K_p(p)/K_p)\cong \Bbb Z_p(1)$ and otherwise we have
$T(K_p(p)/K_p)=1$. The following corollary is driven from this
theorem in 10.5.3 of [Neu-Sch-Wij]

\begin{cor}
With assumptions of above theorem, there is a canonical
isomorphism
$$
 *_{p\in S\setminus T(K_T)} T_p(p)\to Gal(K_S/K_T)(p)
$$
where $T_p$ is the inertia group in the full local group
$G_{K_p}$.
\end{cor}

In the special case of $K=\Bbb Q$, by Kronecker-Weber theory, we
have $\Bbb Q_{S_p\cup S_{\infty}}(p)=\Bbb Q_{\infty}$. Applying
Riemann's existence theorem, one gets the following result of
Neukirch (theorem 10.5.6 in [Neu-Sch-Wij]):

\begin{thm}
Let $p$ be an odd prime, and let $S\supset S_p\cup S_{\infty}$ be
a finite set of places of $\Bbb Q$, then the canonical
homomorphism
$$
 *_{p\in S\setminus S_p\cup S_{\infty}} T(\Bbb Q_p(p)/\Bbb Q_p)
 \to Gal(\Bbb Q_S(p)/\Bbb Q_{\infty})
$$
is an isomorphism.
\end{thm}
Note that $Gal(\Bbb Q_{\infty}/\Bbb Q)\cong \Bbb Z_p$. \\ $\bold
3.$ The following free product decomposition theorem is also
proven 10.7.2 of [Neu-Sch-Wij].

\begin{thm}
Let $\delta=1$ (or $\delta_v=1$), if $\mu_p\subset K$ (or
$\mu_p\subset K_v$ respectively). Suppose $\delta=\sum_v
\delta_v$ and that for $S_0\subseteq S^{finite}$ the dual group of
$$
ker(H^1(Gal_S(K),\mu_p)\to \bigoplus_{v\in S_0} H^1(K_v,\mu_p))
$$
vanishes, then there exists a set of primes $T\supseteq S$ such
that the canonical homomorphism
$$
*_{v\in S\setminus S_0} Gal(K_v(p)/K_v) \ast *_{v\in S\setminus
S_0} Gal(K_v(p)/K_v)/T_v \to Gal_S(K)(p)
$$
is an isomorphism, where $T_v$ denotes the Inertia subgroup of
$v$.
\end{thm}

For $K=\Bbb Q$, $p$ arbitrary or for $K=\Bbb Q(\mu_p)$, $p$
regular with $S= S_p\cup S_{\infty}$, the assumptions on the
kernel above are satisfied, and we get another freeness statement.

\section{Diophantine results}

It is time to prove our main result.
\\ \textbf{Proof(Main Theorem).} By Shafarevich's
conjecture $Gal(k^{cycl}/k)$ is a free profinite group. Let us
call it $F$. The following exact sequence
$$
0 \To \pi_1(X_{\bar {k}}) \To \pi_1(X_{k^{cycl}}) \To Gal(\bar {k}/k^{cycl}) \To 0
$$
is uniquely associated to an outer representation
$$
Gal(\bar {k}/k^{cycl})\to Out(\pi_1(X_{\bar {k}})).
$$
This is because $\pi_1(X_{\bar {k}})$ has trivial center(Look at
[Eil-Mac] and chapter IV of [Bro]). The above exact sequence embeds in the following
$$
0 \To \pi_1(X_{\bar {k}}) \To \pi_1(X_k) \To Gal(\bar {k}/k) \To 0
$$
the quotient being $F$. Now, any point over $k$ gives a spliting of the above short exact sequence,
which in turn induces an splitting of the upper short exact sequence.
Conjugating by inner automorphisms of $Gal(\bar {k}/k)$ which lift inner automorphisms of $F$,
we can induce infinitely many splittings of
$$
0 \To \pi_1(X_{\bar {k}}) \To \pi_1(X_{k^{cycl}}) \To Gal(\bar {k}/k^{cycl}) \To 0
$$
which by section conjecture correspont to points defined over $k^{cycl}$.$\Box$

\begin{rem}
As mentioned before, there exists a pro-$p$ version of the Grothendieck short exact
sequence for a curve defined over any finitely generated field
$$
0 \To \pi^{(p)}_1(X_{\bar k}) \To \pi^{(p)}_1(X_k) \To
Gal(\bar{k}/k)(p) \To 0
$$
and the philosophy of Grothendieck's section conjecture could be
also implied in this case. In the previous section, many examples
are introduced where we get freeness of certain pro-$p$ Galois
groups for $p$ outside a finite set of primes $S$. Together with
hyperbolicity of the corresponding $\pi^{(p)}_1(X_{\bar k})$ one
is lead to use Grothendieck's fudamental conjecture to prove
Diophantine results out of a pro-$p$ version of the Grothendieck's
conjectures look at [Pop2].
\end{rem}
\section*{Acknowledgements}

I would like to thank Abdus Salam International Center for
Theoritical physics (ICTP) and Institute des hautes etude scientifique (IHES) for hospitality and D.T. Le for constant
encouragement in the process of preparation of this note. I would
also like to thank F. Oort and M. Gromov for conversations which
led indirectly to this result. I shall thank corrections to an earlier version suggested by M. Hadian and N. Nikolov.


\end{document}